\documentclass[12pt]{article}
\input{psfig.sty}

\begin{document}
\title{Machines, Logic and Quantum Physics}
\author{David Deutsch and Artur Ekert\\
{\small \emph{Centre for Quantum Computation, Clarendon Laboratory, }}\\
{\small \emph{University of Oxford, Parks Road, Oxford, OX1 3PU, U.K.}}
\and Rossella Lupacchini\\
{\small \emph{Dipartimento di Filosofia, Universita di Bologna, }}\\
{\small \emph{Via Zamboni 38, 40126 Bologna, Italy.}}}
\date{19 November 1999}
\maketitle

\begin{abstract}
Though the truths of logic and pure mathematics are objective and
independent of any contingent facts or laws of nature, our
{\em knowledge} of these truths depends entirely on our knowledge of
the laws of physics. Recent progress in the quantum theory of
computation has provided practical instances of this, and forces us
to abandon the classical view that computation, and hence
mathematical proof, are purely logical notions independent of that
of computation as a physical process. Henceforward, a proof must be
regarded not as an abstract object or process but as a physical
process, a species of computation, whose scope and reliability
depend on our knowledge of the physics of the computer concerned.
\end{abstract}

\section{Mathematics and the physical world}

Genuine scientific knowledge cannot be certain, nor can it be
justified a priori. Instead, it must be conjectured, and then
tested by experiment, and this requires it to be expressed in a
language appropriate for making precise, empirically testable
predictions. That language is mathematics.

This in turn constitutes a statement about what the physical world
must be like if science, thus conceived, is to be possible. As
Galileo put it, ``the universe is written in the language of
mathematics''\cite{G623}. Galileo's introduction of mathematically
formulated, testable theories into physics marked the transition
from the Aristotelian conception of physics, resting on supposedly
necessary a priori principles, to its modern status as a
theoretical, conjectural and empirical science. Instead of seeking
an infallible universal mathematical design, Galilean science uses
mathematics to express quantitative descriptions of an objective
physical reality. Thus mathematics has become the language in which
we express our knowledge of the physical world. This language is
not only extraordinarily powerful and precise, but also effective
in practice. Eugene Wigner referred to ``the unreasonable
effectiveness of mathematics in the physical sciences''\cite{W60}.
But is this effectiveness really unreasonable or miraculous?

Consider how we learn about mathematics. Do we -- that is, do our
brains -- have direct access to the world of abstract concepts and
the relations between them (as Plato believed, and as Roger Penrose
now advocates\cite{PE94}), or do we learn mathematics by
experience, that is by interacting with physical objects? We
believe the latter. This is not to say that the subject-matter of
mathematical theories is in any sense part of, or emerges out of,
the physical world. We do not deny that numbers, sets, groups and
algebras have an autonomous reality quite independent of what the
laws of physics decree, and the properties of these mathematical
structures are just as objective as Plato believed they were. But
they are revealed to us only through the physical world. It is only
physical objects, such as computers or human brains, that ever give
us glimpses of the abstract world of mathematics. But how?

It is a familiar fact, and has been since the prehistoric
beginnings of mathematics, that simple physical systems like
fingers, tally sticks and the abacus can be used to represent some
mathematical entities and operations. Historically the operations
of elementary arithmetic were also the first to be delegated to
more complex machines. As soon as it became clear that additions
and multiplications can be performed by a sequence of basic
procedures and that these procedures are implementable as physical
operations, mechanical devices designed by Blaise Pascal, Gottfried
Wilhelm Leibniz and others began to relieve humans from tedious
tasks such as multiplying two large integers \cite{G72}. In the
twentieth century, following this conquest of arithmetic, the
logical concept of computability was the next to be delegated to
machines. Turing machines were invented in order to formalise the
notion of ``effectiveness'' inherent in the intuitive idea of
calculability. Alan Turing conjectured that the theoretical
machines in terms of which he defined computation are capable of
performing any \emph{finite}, \emph{effective} procedure
(algorithm). It is worth noting that Turing machines were intended
to reproduce every definite operation that a \emph{human} computer
could perform, following preassigned instructions. Turing's method
was to think in terms of physical operations, and imagine that
every operation performed by the computer ``consists of some change
of the physical system consisting of the computer and his
tape''\cite{T36}. The key point is that since the outcome is not
affected by constructing ``a machine to do the work of this
computer'', the effectiveness of a human computer can be mimicked
by a logical machine.

The Turing machine was an abstract construct, but thanks to
subsequent developments in the theory of computation, algorithms
can now be performed by real automatic computing machines. The
natural question now arises: what, precisely, is the set of logical
procedures that can be performed by a physical device? The theory
of Turing machines cannot, even in principle, answer this question,
nor can any approach based on formalising traditional notions of
effective procedures. What we need instead is to extend Turing's
idea of \emph{mechanising} procedures, in particular, the
procedures involved in the notion of derivability. This would
define mathematical proofs as being mechanically reproducible and
to that extent effectively verifiable. The universality and
reliability of logical procedures would be guaranteed by the
mechanical procedures that effectively perform logical operations
-- but by no more than that. But what does it mean to involve real,
physical machines in the definition of a logical notion? and what
might this imply \emph{in return} about the `reasonableness' or
otherwise of the effectiveness of physics in the mathematical
sciences?

While the abstract model of a machine, as used in the classical
theory of computation, is a pure-mathematical construct to which we
can attribute any consistent properties we may find convenient or
pleasing, a consideration of actual computing machines as physical
objects must take account of their actual physical properties, and
therefore, in particular, of the laws of physics. Turing's machines
(with arbitrarily long tapes) can be built, but no one would ever
do so except for fun, as they would be extremely slow and
cumbersome. We find the computers now available much faster and
more reliable. Where does this reliability come from? How do we
know that the computer generates the same outputs as the
appropriate abstract Turing machine, that the machinery of
cog-wheels must finally display the right answer? After all, nobody
has tested the machine by following all possible logical steps, or
by performing all the arithmetic it can perform. If they were able
and willing to do that, there would be no need to build the
computer in the first place. The reason we trust the machine cannot
be based entirely on logic; it must also involve our knowledge of
the physics of the machine. We take for granted the laws of physics
that govern the computation, i.e. the physical process that takes
the machine from an initial state (input) to a final state
(output). Moreover, our understanding is informed by physical
theories which, though formulated in mathematical terms in the
tradition of Galileo, evolved by conjectures and empirical
refutations. In this perspective, what Turing really asserted was
that it is possible to build a universal computer, a machine that
can be programmed to perform any computation that any other
physical object can perform. That is to say, a single buildable
physical object, given only proper maintenance and a supply of
energy and additional memory when needed, can mimic all the
behavior and responses of any other physically possible object or
process. In this form Turing's conjecture (which Deutsch has called
in this context the \emph{Church-Turing principle}~\cite{D85}) can
be viewed as a statement about the physical world.

Are there any limits to computations performed by computing
machines? Obviously there are both logical and physical limits.
Logic tells us that, for example, no machine can find more than one
even prime, whilst physics tells us that, for example, no
computations can violate the laws of thermodynamics. Moreover,
logical and physical limitations can be intimately linked, as
illustrated by the ''halting problem'' . According to logic, the
halting problem says that there is no algorithm for deciding
whether any given machine, when started from any given initial
situation, eventually stops. Therefore some computational problems,
such as determining whether a specified universal Turing machine,
given a specified input, will halt, cannot be solved by any Turing
machine. In physical terms, this statement says that machines with
certain properties cannot be physically built, and as such can be
viewed as a statement about physical reality or equivalently about
the laws of physics.

So where does mathematical effectiveness come from? It is not
simply a miracle, ``a wonderful gift which we neither understand
nor deserve''~\cite{W60} -- at least, no more so than our ability
to discover \emph{empirical} knowledge, for our knowledge of
mathematics and logic is inextricably entangled with our knowledge
of physical reality: every mathematical proof depends for its
acceptance upon our agreement about the rules that govern the
behavior of physical objects such as computers or our brains. Hence
when we improve our knowledge about physical reality, we may also
gain new means of improving our knowledge of logic, mathematics and
formal notions. It seems that we have no choice but to recognize
the dependence of our mathematical \emph{knowledge} (though not, we
stress, of mathematical truth itself) on physics, and that being
so, it is time to abandon the classical view of computation as a
purely logical notion independent of that of computation as a
physical process. In the following we discuss how the discovery of
quantum mechanics in particular has changed our understanding of
the nature of computation.

\section{Quantum interference}

To explain what makes quantum computers so different from their classical
counterparts, we begin with the phenomenon of quantum interference. Consider
the following computing machine whose input can be prepared in one of two
states representing, $0$ and $1$.

\begin{figure}[ht]
\vskip .1cm
\centerline{
\psfig{file=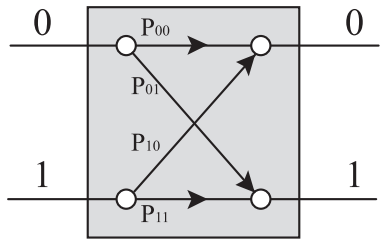,width=50mm}}
\vskip .2cm
\caption{Schematic representation of the most general machine that performs
a computation mapping $\left\{ 0,1\right\}$ to itself. Here $p_{ij}
$ is the probability for the machine to produce the output $j$ when
presented with the input $i$. (The action of the machine depends on
no other input or stored information.)}
\label{qmachine}
\end{figure}

The machine has the property that if we prepare its input with the
value $a$ ($a=0$ or $1$) and then measure the output, we obtain,
with probability $p_{ab} $, the value $b$ ($b=0$ or $1$). It may
seem obvious that if the $p_{ab}$ are arbitrary apart from
satisfying the standard probability conditions $\sum_{b}p_{ab}=1$,
Fig. (\ref{qmachine}) represents \emph{the most general} machine
whose action depends on no other input or stored information and
which performs a computation mapping $\left\{ 0,1\right\} $ to
itself. The two possible deterministic limits are obtained by
setting $p_{01}=$ $p_{10}=0,\ p_{00}=$ $p_{11}=1$ (which gives a
logical identity machine) or $p_{01}=$ $p_{10}=1,$ $p_{00}=$
$p_{11}=0$ (which gives a negation (`{\bf not}') machine).
Otherwise we have a randomising device. Let us assume, for the sake
of illustration, that $p_{01}=$ $p_{10}=p_{00}=$ $p_{11}=0.5.$
Again, we may be tempted to think of such a machine as a random
switch which, with equal probability, transforms any input into one
of the two possible outputs. However, that is not necessarily the
case. When the particular machine we are thinking of is followed by
another, identical, machine the output is always the negation of
the input.

\begin{figure}[ht]
\vskip .1cm
\centerline{
\psfig{file=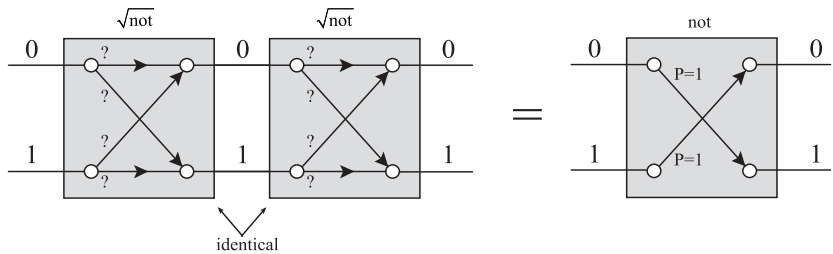,width=110mm}}
\vskip .2cm
\caption{Concatenation of the two identical machines mapping $\left\{
0,1\right\} $ to itself. Each machine, when tested separately,
behaves as a random switch, however, when the two machines are
concatenated the randomness disappears - the net effect is the
logical operation {\bf not}. This is in clear contradiction with
the axiom of additivity in probability theory! }
\label{2qmachines}
\end{figure}

This is a very counter-intuitive claim - the machine alone outputs
$0$ or $1$ with equal probability and independently of the input,
but the two machines, one after another, acting independently,
implement the logical operation {\bf not}. That is why we call this
machine $\sqrt{{\rm \bf not}}.$ It may seem reasonable to argue
that since there is no such operation in logic, the $\sqrt{{\rm \bf
not}}$ machine cannot exist. But it does exist! Physicists studying
single-particle interference routinely construct them, and some of
them are as simple as a half-silvered mirror i.e. a mirror which
with probability 50\% reflects a photon which impinges upon it and
with probability 50\% allows it to pass through. Thus the two
concatenated machines are realised as a sequence of two
half-silvered mirrors with a photon in each denoting $0$ if it is
on one of the two possible paths and $1$ if it is on the other.

\begin{figure}[ht]
\vskip .1cm
\centerline{
\psfig{file=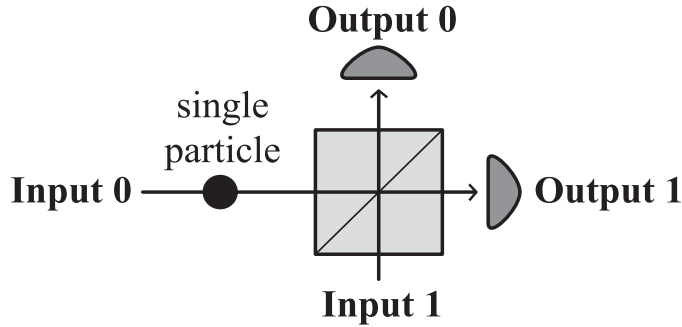,width=90mm}}
\vskip .2cm
\caption{The experimental realisation of the $\sqrt{\rm\bf not}$ gate. A
half-silvered mirror reflects half the light that impinges upon it.
But a single photon doesn't split: when we send a photon at such a
mirror it is detected, with equal probability, either at Output 0
or 1. This does not, however, mean that the photon leaves the
mirror in either of the two outputs at random. In fact the photon
takes both paths at once! This can be demonstrated by concatenating
two half-silvered mirrors as shown in the next figure.}
\label{exper1}
\end{figure}

\begin{figure}[ht]
\vskip .1cm
\centerline{
\psfig{file=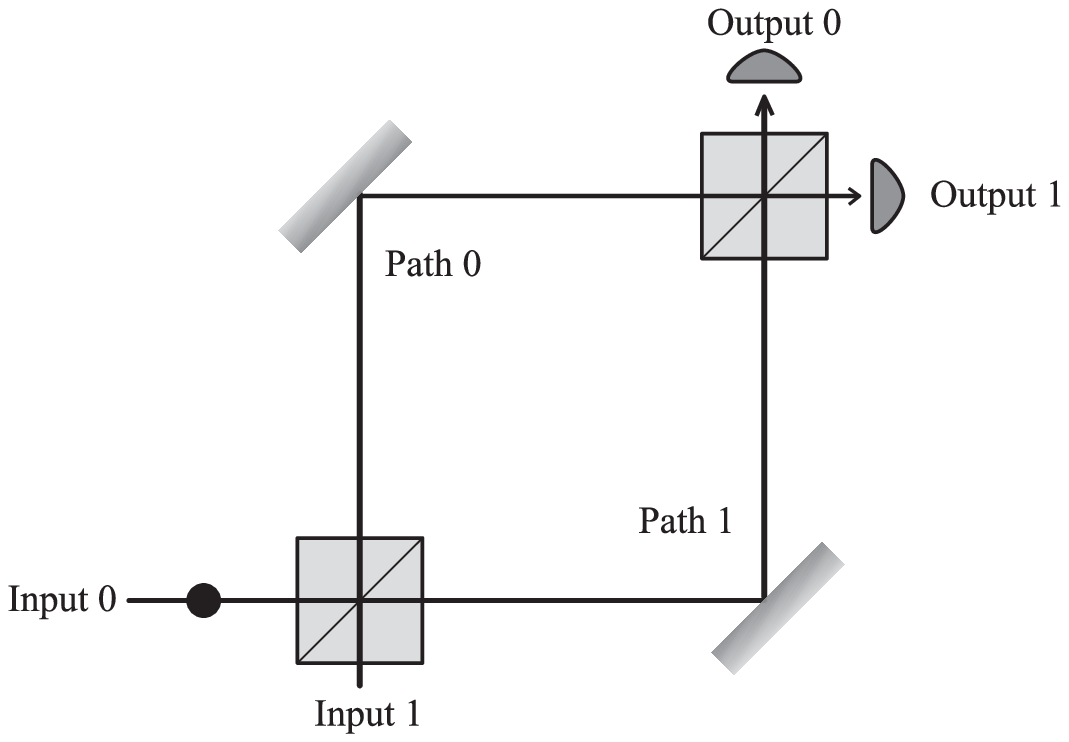,width=110mm}}
\vskip .2cm
\caption{The experimental realisation of the two concatenated $\sqrt{\rm\bf
not}$ gates, known as a single-particle interference. A photon
which enters the interferometer via Input 0 always strikes a
detector Output 1 and never a detector at Output 0. Any explanation
which assumes that the photon takes exactly one path through the
interferometer leads to the conclusion that the two detectors
should on average each fire on half the occasions when the
experiment is performed. But experiment shows otherwise!}
\label{exper2}
\end{figure}

The reader may be wondering what has happened to the axiom of
additivity in probability theory, which says that if $E_{1}$ and
$E_{2}$ are mutually exclusive events then the probability of the
event $(E_{1}$ \textbf{or} $E_{2})$ is the sum of the probabilities
of the constituent events, $E_{1}$, $E_{2}$. We may argue that the
transition $0\rightarrow 0$ in the composite machine can happen in
the two mutually exclusive ways, namely, $0\rightarrow 0\rightarrow
0$ or $0\rightarrow 1\rightarrow 0.$ The probabilities of the two
are $p_{00}p_{00}$ and $p_{01}p_{10}$ respectively. Thus the sum
$p_{00}p_{00}+p_{01}p_{10}$ represents the probability of the
$0\rightarrow 0$ transition in the new machine. Provided that
$p_{00}$ or $p_{01}p_{10}$ are different from zero, this
probability should also be different from zero. Yet we can build
machines in which $p_{00}$ and $p_{01}p_{10}$ are different from
zero, but the probability of the $0\rightarrow 0$ transition in the
composite machine is equal to zero. So what is wrong with the above
argument?

One thing that is wrong is the assumption that the processes
$0\rightarrow 0\rightarrow 0$ and $0\rightarrow 1\rightarrow 0$ are
mutually exclusive. In reality, the two transitions both occur,
simultaneously. We cannot learn about this fact from probability
theory or any other a priori mathematical construct. We learn it
from the best physical theory available at present, namely quantum
mechanics. Quantum theory explains the behavior of $\sqrt{{\rm\bf
not}}$ and correctly predicts the probabilities of all the possible
outputs no matter how we concatenate the machines. This knowledge
was created as the result of conjectures, experimentation, and
refutations. Hence, reassured by the physical experiments that
corroborate this theory, logicians are now entitled to propose a
new logical operation $\sqrt{{\rm\bf not}}$. Why? Because a
faithful physical model for it exists in nature!

Let us now introduce some of the mathematical machinery of quantum
mechanics which can be used to describe quantum computing machines
ranging from the simplest, such as $\sqrt{{\rm\bf not}}$ , to the
quantum generalisation of the universal Turing machine. At the
level of predictions, quantum mechanics introduces the concept of
\emph{probability amplitudes} -- complex numbers $c $ such that the
quantities $\left| c\right| ^{2}$ may under suitable circumstances
be interpreted as probabilities. When a transition, such as ``a
machine composed of two identical sub-machines starts in
state $0$ and generates output $0$, and nothing else happens'', can
occur in several alternative ways, the overall probability
amplitude for the transition is the sum, not of the probabilities,
but of the probability amplitudes for each of the constituent
transitions considered separately.

\begin{figure}[ht]
\vskip .1cm
\centerline{
\psfig{file=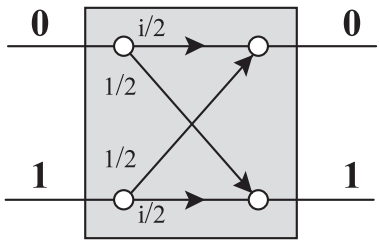,width=50mm}}
\vskip .2cm
\caption{Transitions in quantum machines are described by probability
amplitudes rather than probabilities. Probability amplitudes are
complex numbers $c$ such that the quantities $\left| c\right| ^{2}$
may under suitable circumstances be interpreted as probabilities.
When a transition, such as ``a machine composed of two
identical sub-machines starts in state $0$ and generates output
$0$, and nothing else happens'', can occur in several alternative
ways, the probability amplitude for the transition is the sum of
the probability amplitudes for each of the constituent transitions
considered separately.}
\label{interference}
\end{figure}

In the $\sqrt{{\rm\bf not}}$ machine, the probability amplitudes of
the $0\rightarrow 0$ and $1\rightarrow 1$ transitions are both
$i/\sqrt{2}$, and the probability amplitudes of the $0\rightarrow
1$ and $1\rightarrow 0$ \noindent transitions are both
$1/\sqrt{2}$. This means, for example, that the $\sqrt{{\rm\bf
not}}$ machine preserves the bit value with the probability
amplitude $c_{00}=c_{11}=i/\sqrt{2}$ and negates it with the
probability amplitude $c_{01}=c_{10}=1/\sqrt{2}.$ In order to
obtain the corresponding probabilities we have to take the modulus
squared of the probability amplitudes which gives $1/2$ both for
preserving and swapping the bit value. This describes the behavior
of the $\sqrt{{\rm\bf not}}$ machine in Fig.(\ref{qmachine}).
However, when we concatenate the two machines, as in
Fig.(\ref{2qmachines}) then, in order to calculate the probability
of output $0$ on input $0$, we have to add the probability
amplitudes of all computational paths leading from input $0$ to
output $0.$ There are only two of them - $c_{00}c_{00}$ and
$c_{01}c_{10}$. The first computational path has probability
amplitude $i/\sqrt{2}\times i/\sqrt{2}=-1/2$ and the second one
$1/\sqrt{2}\times1/\sqrt{2}=+1/2.$ We add the two probability
amplitudes first and then we take the modulus squared of the sum.
We find that the probability of output $0$ is zero. Unlike
probabilities, probability amplitudes can cancel each other out!

\section{Quantum algorithms}

Addition of probability amplitudes, rather then probabilities, is
one of the fundamental rules for prediction in quantum mechanics
and applies to all physical objects, in particular quantum
computing machines. If a computing machine starts in a specific
initial configuration (input) then the probability that after its
evolution via a sequence of intermediate configurations it ends up
in a specific final configuration (output)is the squared modulus of
the sum of all the probability amplitudes of the computational
paths that connect the input with the output. The amplitudes are
complex numbers and may cancel each other, which is referred to as
\emph{destructive interference}, or enhance each other, referred to
as \emph{constructive interference}. The basic idea of quantum
computation is to use quantum interference to amplify the correct
outcomes and to suppress the incorrect outcomes of computations.
Let us illustrate this by describing a variant of the first quantum
algorithm, proposed by David Deutsch in 1985.

\begin{figure}[ht]
\vskip .1cm
\centerline{
\psfig{file=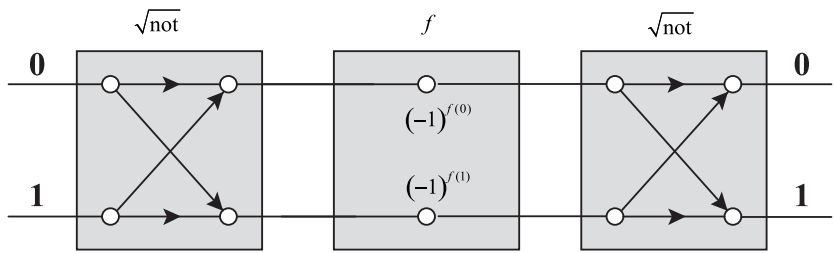,width=110mm}}
\vskip .2cm
\caption{Schematic representation of a quantum machine that solves Deutsch's
problem with a single function evaluation.}
\label{function}
\end{figure}

Consider the Boolean functions $f$ that map $\{0,1\}$ to $\{0,1\}$.
There are exactly four such functions: two constant functions
($f(0)=f(1)=0$ and $f(0)=f(1)=1$) and two ``balanced'' functions
($f(0)=0,f(1)=1$ and $f(0)=1,f(1)=0$). Suppose we are allowed to
evaluate the function $f$ \emph{only once} (given, say, a lengthy
algorithm for evaluating it on a given input, or a look-up table
that may be consulted only once) and asked to determine whether $f$
is constant or balanced (in other words, whether $f(0)$ and $f(1)$
are the same or different). Note that we are not asking for the
particular values $f(0)$ and $f(1)$ but for a global property of
the function $f$. Our classical intuition insists, and the
classical theory of computation confirms, that to determine this
global property of $f$, we have to evaluate both $f(0)$ and $f(1)$,
which involves evaluating $f$ twice. Yet this is simply not true in
physical reality, where quantum computation can solve Deutsch's
problem with a single function evaluation. The machine that solves
the problem, using quantum interference, is composed of the two
$\sqrt{{\rm\bf not}}$s with the function evaluation machine in
between them, as in Fig.(\ref{function}).

We need not go into the technicalities of the physical
implementation of the evaluation of $f(x)$, where $f$ is in general
a Boolean function mapping $\left\{0,1\right\}^{n}\rightarrow
\left\{ 0,1\right\} ^{m}$. But generally, a machine that evaluates
such a function must be capable of traversing as many computational
paths as there are possible values $x$ in the domain of $f$ (so we
can label them with $x$). Its effect is that the probability
amplitude on path $x$ is multiplied by the phase factor $\exp\left(
\frac{2\pi i f(x)}{2^{m}}\right) $ \cite{CEMM98}. In the case of
Deutsch's problem the two phase factors are $\left(-1\right)^{
f(0)}$ and $\left( -1\right)^{f(1)}$. Now we can calculate the
probability amplitude of output $0$ on input $0.$ The probability
amplitudes on the two different computational paths are $i/\sqrt{2}
\times\left( -1\right)^{f(0)} \times i/\sqrt{2}=-1/2
\times(-1)^{f(0)}$ and $1/\sqrt{2}\times\left( -1\right)
^{f(1)}\times1/\sqrt{2}=1/2\times(-1)^{f(1)}.$ Their sum is

\begin{equation}
\frac{1}{2}\left(  (-1)^{f(1)}-(-1)^{f(0)}\right),
\end{equation}

which is $0$ when $f$ is constant and $\pm1$ when $f$ is balanced.
Thus the probability of output $0$ on input $0$ is given by the
modulus squared of the expression above, which is zero when $f$ is
constant and unity when $f$ is balanced. Deutsch's result laid the
foundation for the new field of quantum computation. The hunt began
for useful tasks for quantum computers to do. A sequence of
steadily improved quantum algorithms led in 1994 to Peter Shor's
discovery of a quantum algorithm that, in principle, could perform
efficient factorisation~\cite{S94}. Since the intractability of
factorisation underpins the security of many of the most secure
known methods of encryption, including the most popular public key
cryptosystem RSA~\cite{RSA}~\footnote{In December 1997 the British
Government officially confirmed that public-key cryptography was
originally invented at the Government Communications Headquarters
(GCHQ) in Cheltenham. By 1975, James Ellis, Clifford Cocks, and
Malcolm Williamson from GCHQ had discovered what was later
re-discovered in academia and became known as RSA and
Diffie-Hellman key exchange.}, Shor's algorithm was soon hailed as
the first `killer application' for quantum computation
--- something very useful that only a quantum computer could do.

Few if any mathematicians doubt that the factoring problem is in
the ``\textbf{BPP''}class (where \textbf{BPP} stands for ``bounded
error probabilistic polynomial time''), but interestingly, this has
never been proved\textbf{. }\noindent In computational complexity
theory it is customary to view problems in \textbf{BPP} as being
``tractable'' or ``solvable in practice'' and problems not in
\textbf{BPP} as ``intractable'' or ``unsolvable in practice on a
computer'' (see, for example, \cite{PA94}). A `\textbf{BPP
}\noindent class algorithm' for solving a problem is an efficient
algorithm which, for any input, provides an answer that is correct
with a probability greater than some constant $\delta>1/2$. In
general we cannot check easily if the answer is correct or not but
we may repeat the computation some fixed number $k$ times and then
take a majority vote of all the $k$ answers. For sufficiently large
$k$ the majority answer will be correct with probability as close
to 1 as desired. Now, Shor's result proves that factoring is not in
reality an intractable task -- and we learned this by studying
quantum mechanics!

As a matter of fact, Richard Feynman, in his talk during the First
Conference on the Physics of Computation held at MIT in 1981,
observed that it appears to be impossible to simulate a general
quantum evolution on a classical probabilistic computer in an
efficient way~\cite{FY82}. That is to say, any classical simulation
of quantum evolution involves an exponential slowdown in time
compared with the natural evolution, since the amount of
information required to describe the evolving quantum state in
classical terms generally grows exponentially with time. However,
instead of viewing this fact as an obstacle, Feynman regarded it as
an opportunity. If it requires so much computation to work out what
will happen in a complicated multiparticle interference experiment
then, he argued, the very act of setting up such an experiment and
measuring the outcome is tantamount to performing a complex
computation. Thus, Feynman suggested that it might be possible to
simulate a quantum evolution efficiently after all, provided that
the simulator itself is a quantum mechanical device. Furthermore,
he conjectured that if one wanted to simulate a different quantum
evolution, it would not be necessary to construct a new simulator
from scratch. It should be possible to choose the simulator so that
minor systematic modifications of it would suffice to give it any
desired interference properties. He called such a device a
universal quantum simulator. In 1985 Deutsch proved that such a
universal simulator or a universal quantum computer does exist and
that it could perform any computation that any other quantum
computer (or any Turing-type computer) could perform \cite{D85}.
Moreover, it has since been shown that the time and other resources
that it would need to do these things would not increase
exponentially with the size or detail of the physical system being
simulated, so the relevant computations would be tractable by the
standards of complexity theory~\cite{BV93}. This illustrates the
fact that the more we know about physics, the more we know about
computation and mathematics. Quantum mechanics proves that
factoring is tractable: without quantum mechanics we do not yet
know how to settle the issue either way.

\section{Deterministic, Probabilistic, and Quantum Computers}

Any quantum computer, including the universal one, can be described
in a fashion similar to the special-purpose machines we have
described above, essentially by replacing probabilities by
probability amplitudes. Let us start with a classical Turing
machine. This is defined by a finite set of quintuples of the form
\begin{equation}
(q,s,q^{\prime},s^{\prime},d), \label{TM}
\end{equation}
where the first two characters describe the initial condition at
the beginning of a computational step and the remaining three
characters describe the effect of the instruction to be executed in
that condition ($q$ is the current configuration, $s$ is the symbol
currently scanned, $q^{\prime}$ is the configuration to enter next,
$s^{\prime}$ is the symbol to replace $s$, and $d$ indicates motion
of one square to the right, or one square to the left, or stay
fixed, relative to the tape). In this language a computation
consists of presenting the machine with an input which is a finite
string of symbols from the alphabet $\Sigma$ written in the tape
cells, then allowing the machine to start in the initial state
$q_{0}$ with the head scanning the leftmost symbol of the input and
to proceed with its basic operations until it stops in its final
(halting) state $q_{h}$.(In some cases the computation might not
terminate.) The output of the computation is defined as the
contents of some chosen part of the tape when (and if) the machine
reaches its halting state.

During a computation the machine goes through a sequence of
configurations; each configuration provides a global description of
the machine and is determined by the string written on the entire
tape, the state of the head and the position of the head. For
example, the initial configuration is given by the input string,
state $q_{0}$, and the head scanning the leftmost symbol from the
input. There are infinitely many possible configurations of the
machine but in all successful computations the machine goes through
only a finite sequence of them. The transitions between
configurations are completely determined by the quintuples
(\ref{TM}).

\begin{figure}[h]
\vskip .1cm
\centerline{
\psfig{file=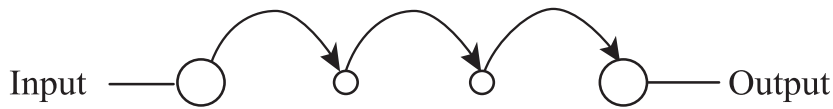,width=110mm}}
\vskip .2cm
\caption{A three step deterministic computation.}
\label{deterministic}
\end{figure}

Computations do not, in principle, have to be deterministic.
Indeed, we can augment a Turing machine by allowing it ``to toss an
unbiased coin'' and to choose its steps randomly. Such a
probabilistic computation can be viewed as a directed, tree-like
graph where each node corresponds to a configuration of the
machine, and each edge represents one step of the computation. The
computation starts from the root node representing the initial
configuration and it subsequently branches into other nodes
representing configurations reachable with non-zero probability
from the initial configuration. The action of the machine is
completely specified by a finite description of the form:

\begin{equation}
\delta:\quad Q\times\Sigma\times Q\times\Sigma\times
\{{\rm Left,Right,Nothing}\}\longmapsto\lbrack0,1],
\label{PTM}
\end{equation}
where $\delta(q,s,q^{\prime},s^{\prime},d)$ gives the probability
that if the machine is in state $q$ reading symbol $s$ it will
enter state $q^{\prime}$ , write $s^{\prime}$ and move in direction
$d$. This description must conform to the laws of probability as
applied to the computation tree. If we associate with each edge of
the graph the probability that the computation follows that edge
then we must require that the sum of the probabilities on edges
leaving any single node is always equal to $1$. Probability of a
particular path being followed from the root to a given node is the
product of the probabilities along the path's edges, and the
probability of a particular configuration being reached after $n$
steps is equal to the sum of the probabilities along all paths
which in $n$ steps connect the initial configuration with that
particular configuration. Such randomized algorithms can solve some
problems (with arbitrarily high probability less than 1) much
faster than any known deterministic algorithms.

\begin{figure}[ht]
\vskip .1cm
\centerline{
\psfig{file=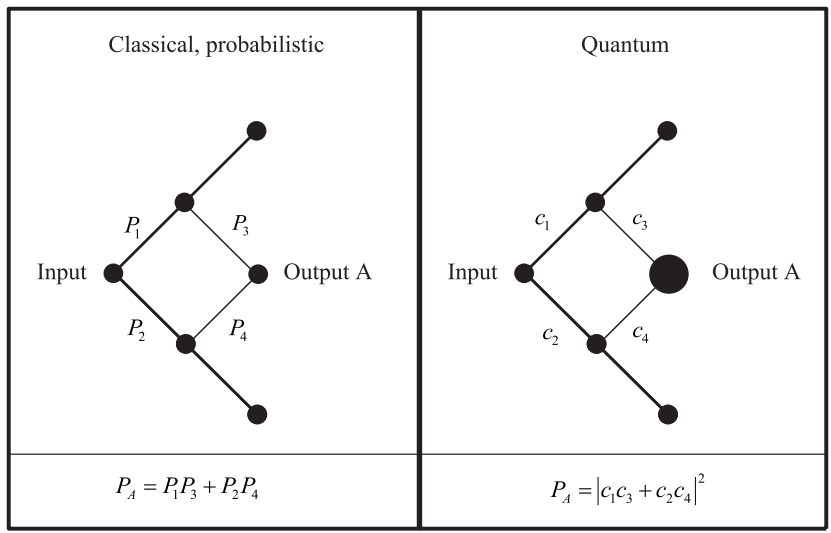,width=110mm}}
\vskip .2cm
\caption{The probabilistic Turing machine (left) - the probability of output
A is the sum of the probabilities of all computations leading to
output A. In the quantum Turing machine (on the right) the
probability of output A is obtained by adding all probability
amplitudes leading from the initial state to output A and then
taking the squared modulus of the sum. In the quantum case
probabilities of some outcomes can be enhanced (constructive
interference) or suppressed (destructive interference) compared
with what classical probability theory would permit.}
\label{PTMQTM}
\end{figure}

The classical model described above suggests a natural quantum
generalisation. A quantum computation can be represented by a graph
similar to that of a probabilistic computation. Following the rules
of quantum dynamics we associate with each edge in the graph the
probability \emph{amplitude} that the computation follows that
edge. As before, the probability amplitude of a particular path
being followed is the product of the probability amplitudes along
the path's edges and the probability amplitude of a particular
configuration is the sum of the amplitudes along all possible paths
leading from the root to that configuration. If a particular final
configuration can be reached via two different paths with
amplitudes $c$ and $-c$ then the probability of reaching that
configuration is $|c-c|^{2}=0$ despite the fact that the
probability for the computation to follow either of the two paths
separately is $|c|^{2}$ in both cases. Furthermore a single quantum
computer can follow many distinct computational paths
simultaneously and produce a final output depending on the
interference of all of them. This is in contrast to a classical
probabilistic Turing machine which follows only some \emph{single}
(randomly chosen) path. The action of any such quantum machine is
completely specified by a finite description of the form

\begin{equation}
\delta:\quad Q\times\Sigma\times Q\times\Sigma\times
\{{\rm Left,Right,Nothing}\}\longmapsto\mathbf{C} \label{QTM}
\end{equation}
where $\delta(q,s,q^{\prime},s^{\prime},d)$ gives the probability
\emph{amplitude} that if the machine is in state $q$ reading symbol $s$ it
will enter state $q^{\prime}$ , write $s^{\prime}$ and move in direction $d$.

\section{Deeper implications}

When the physics of computation was first investigated, starting in
the 1960s, one of the main motivations was a fear that
quantum-mechanical effects might place fundamental bounds on the
accuracy with which physical objects could render the properties of
the abstract entities, such as logical variables and operations,
that appear in the theory of computation. Thus it was feared that
the power and elegance of that theory, its most significant
concepts such as computational universality, its fundamental
principles such as the Church-Turing thesis and its powerful
results such as the more modern theory of complexity, might all be
mere figments of pure mathematics, not really relevant to anything
in nature.

Those fears turned out to be groundless. Quantum mechanics, far
from placing limits on which Turing computations can be performed
in nature, permits them all, and in addition provides new modes of
computation such as those we have described. As far as the elegance
of the theory goes, it turns out that the quantum theory of
computation hangs together better, and fits in far more naturally
with fundamental theories in other fields, than its classical
approximation was ever expected to. The very word `quantum' means
the same as the word `bit' --- an elementary chunk --- and this
reflects the fact that classical physical systems, being subject to
the generic instability known as `chaos', would not support digital
computation at all (so even Turing machines, the theoretical
prototype of all classical computers, were secretly
quantum-mechanical all along!). The Church-Turing thesis in the
classical theory (that all `natural' models of computation are
essentially equivalent to each other), was never proved. Its
analogue in the quantum theory of computation (the Church-Turing
Principle, that the universal quantum computer can simulate the
behavior of any finite physical system) was straightforwardly
proved in Deutsch's 1985 paper~\cite{D85}. A stronger result (also
conjectured but never proved in the classical case), namely that
such simulations can always be performed in a time that is at most
a polynomial function of the time taken for the physical evolution,
has since been proved in the quantum case \cite{BV93}.

Among the many ramifications of quantum computation for apparently
distant fields of study are its implications for the notion of
mathematical proof. Performing any computation that provides a
definite output is tantamount to \emph{proving} that the observed
output is one of the possible results of the given computation.
Since we can describe the computer's operations mathematically, we
can always translate such a computation into the proof of some
mathematical theorem. This was the case classically too, but in the
absence of interference effects it is always possible to keep a
record of the steps of the computation, and thereby produce (and
check the correctness of) a proof that satisfies the classical
definition - as ``a sequence of propositions each of which is
either an axiom or follows from earlier propositions in the
sequence by the given rules of inference''. Now we are forced to
leave that definition behind. Henceforward, a proof must be
regarded as a process --- the computation itself --- for we must
accept that in future, quantum computers will prove theorems by
methods that neither a human brain nor any other arbiter will ever
be able to check step-by-step, since if the `sequence of
propositions' corresponding to such a proof were printed out, the
paper would fill the observable universe many times over.

\section{Concluding remarks}

This brief discussion has merely scratched the surface of the
rapidly developing field of quantum computation. We have
concentrated on the fundamental issues and have avoided discussing
physical details and technological practicalities. However, it
should be mentioned that quantum computing is a serious possibility
for future generations of computing devices. This is one reason why
the field is now attracting increasing attention from both academic
researchers and industry worldwide. At present it is not clear
when, how and even whether fully-fledged quantum computers will
eventually be built; but notwithstanding this, the quantum theory
of computation already plays a much more fundamental role in the
scheme of things than its classical predecessor did. We believe
that anyone who seeks a fundamental understanding of either
physics, computation or logic must incorporate its new insights
into their world view.

\bigskip
\end{document}